\documentclass[12pt]{article}
\input epsf      
\usepackage{graphicx}
\usepackage{amssymb}
\usepackage{amsfonts}
\usepackage{amsthm}
\usepackage{amsfonts}

\newcommand{\ra}{\rightarrow}

\newcommand{\ba}[1]{\begin{array}{#1}}
\newcommand{\ea}{\end{array}}

\newcommand{\be}{\begin{equation}}
\newcommand{\ee}{\end{equation}}
\newcommand{\bea}{\begin{eqnarray}}
\newcommand{\eea}{\end{eqnarray}}
\newcommand{\beann}{\begin{eqnarray*}}
\newcommand{\eeann}{\end{eqnarray*}}

\def\reff#1{(\ref{#1})}
\newtheorem{conjecture}{Conjecture}

\begin{document}

\bibstyle{ams}

\title{The Length of an SLE - Monte Carlo Studies}


\author{Tom Kennedy
\\Department of Mathematics
\\University of Arizona
\\Tucson, AZ 85721
\\ email: tgk@math.arizona.edu
}
\maketitle

\begin{abstract}

The scaling limits of a variety of critical two-dimensional lattice models
are equal to the Schramm-Loewner evolution (SLE) for a suitable value
of the parameter $\kappa$. These lattice
models have a natural parametrization of their random curves given 
by the length of the curve. This parametrization
(with suitable scaling) should provide a natural parametrization 
for the curves in the scaling limit. We conjecture that this parametrization
is also given by a type of fractal variation along the curve, and 
present Monte Carlo simulations to support this conjecture. Then
we show by simulations that if this fractal variation is used to parametrize
the SLE, then the parametrized curves have the same distribution as 
the curves in the scaling limit of the lattice models with their natural 
parametrization.

\end{abstract}

\bigskip

\newpage

\section{Introduction}
\label{intro}

Schramm-Loewner evolution (SLE) is a one parameter family of random processes
that produce random curves in the plane \cite{schramm}.
(When we refer to the curves of SLE we will mean the trace of 
the SLE.) Beffara \cite{befa,befb} proved that the Hausdorff 
dimension of the SLE curve is $1+\kappa/8$ a.s. for $\kappa \le 8$. 
With the usual definition of length, the length of these SLE curves
is infinite. Nonetheless, one would like to define some notion of 
the length of an SLE curve.  

To motivate our proposal for a definition of length for SLE we first consider
models which come from scaling limits of models on a lattice, e.g., 
the loop-erased random walk (LERW), the self-avoiding walk (SAW), 
interfaces in the critical Ising model, and the percolation exploration 
process. 
Before the scaling limit, the random curves in the lattice model 
(which are walks or interfaces) have a natural 
parametrization arising from the number of steps in the curve, or
equivalently the length of the curve. This leads to a parametrization
of the curves in the scaling limit that we will refer to as the 
``natural parametrization.'' It is defined as follows.
Let $W(n)$ be a random curve in the lattice model with the lattice spacing 
equal to 1. The mean square distance the curve travels after $n$ steps
should grow as $n^{2 p}$ for some critical exponent $p$. 
\be
E \,[|W(n)-W(0)|^2] \approx c n^{2 p}
\ee
The exponent $p$ should be related to the Hausdorff dimension, $d_H$, 
of the curve by $p=1/d_H$. From now on we will assume this is the case
and use $1/d_H$ in place of $p$. 

Initially, $W(n)$ is defined only for non-negative integers $n$. 
We extend $W(t)$  to all positive real $t$ by linearly interpolating.
Then we define
\be
\omega(t) = \lim_{n \rightarrow \infty} n^{-1/d_H} \, W(nt), 
\ee
We will refer to the resulting parametrization of the curves in the 
scaling limit as the natural parametrization.
This definition is analogous to the construction of 
Brownian motion as the scaling limit of the ordinary random walk. 
The goal of this paper is to propose a parametrization of 
SLE which corresponds to the natural parametrization of the 
scaling limits of the various lattice models and to provide some support
from simulations for this correspondence.

To introduce our parametrization of SLE,  
let $0=t_0 < t_1 < t_2 \cdots < t_n = t$ be a partition of the time
interval $[0,t]$. The usual definition of the length or total variation
of the random curve $\gamma$ over the 
time interval $[0,t]$ would be the supremum over all such partitions of 
\be
\sum_{j=1}^n |\gamma(t_j)-\gamma(t_{j-1})|
\ee
The length of $\gamma(t_j)-\gamma(t_{j-1})$ will be of order 
$|t_j-t_{j-1}|^{1/d_H}$, and since $d_H>1$, the total variation of $\gamma$ 
will be infinite. Since the length of a segment is of order 
$|t_j-t_{j-1}|^{1/d_H}$, this suggests that we consider the quantity,
\be
fvar(\gamma[0,t],\Pi)=\sum_{j=1}^n |\gamma(t_j)-\gamma(t_{j-1})|^{d_H}
\ee
where $\Pi$ denotes the partition, $\{t_0,t_1,\cdots,t_n\}$. 
Then we take a sequence of partitions $\Pi_m$ whose mesh converges to zero.
(The mesh of a partition is the length of the largest subinterval.)
Then we consider the limit, 
\be
\lim_{m \ra \infty} fvar(\gamma[0,t],\Pi_m)
\ee
Of course, the existence of this limit and its dependence on the sequence of 
partitions used is a subtle question. 

This quantity is sometimes called the $p$-variation in the stochastic
processes literature. ($p$ is used where we have $d_H$.) 
We find the name $p$-variation particularly dull, and will refer 
to this quantity as the fractal variation. 
When $d_H=2$, this is the quadratic variation 
studied by L\'evy for Brownian motion \cite{levy}. 
It is non-random and proportional to $t$ under suitable conditions 
on the convergence of the sequence of partitions \cite{dud,vega}.

A serious drawback of this definition is that it depends on the 
parametrization of the curve. For example, in the scaling limit of 
a lattice model we could compute this fractal variation using either the 
parametrization of the curve by capacity or its natural parametrization.
Simulations indicate that the two results will not agree. 

\begin{figure}[tbh]
\includegraphics{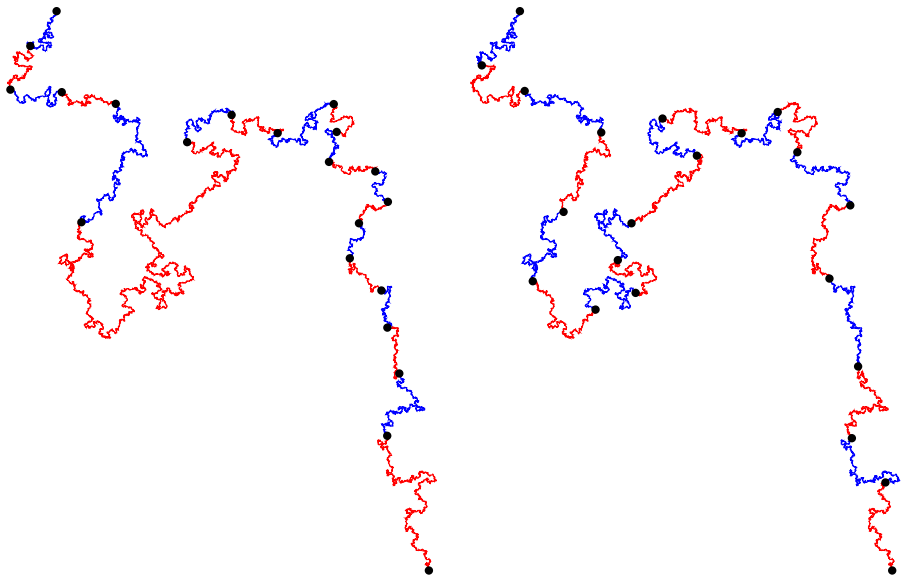}
\caption{\leftskip=25 pt \rightskip= 25 pt 
An SLE$_{8/3}$ curve divided into equal segments according to 
parametrization by capacity (left) and by fractal variation (right).
(Color online.) 
}
\label{sle_segments}
\end{figure}

In our setting the dependence of the above definition on the choice of 
parametrization is particularly troubling since our goal in 
introducing this quantity is to define a parametrization. 
A better definition that does not depend on the choice of parametrization 
of the curve is the following.
Let $\Delta t>0$. We define $t_i$ inductively. Let $t_0=0$. Given
$t_{i-1}$, let $t_i$ be the first time after $t_{i-1}$ with 
\bea
|\gamma(t_i)-\gamma(t_{i-1})|=(\Delta t)^{1/d_H}
\nonumber
\eea
Let $n$ be the last index with $t_n<t$. Then we will define 
\be
fvar(\gamma[0,t],\Delta t)= n \Delta t
= \sum_{j=1}^n |\gamma(t_j)-\gamma(t_{j-1})|^{d_H} 
\ee
and define the fractal variation of $\gamma$ over $[0,t]$ to be
\be
fvar(\gamma[0,t]) = \lim_{\Delta t \ra 0} \, fvar(\gamma[0,t],\Delta t)
\label{fvar}
\ee

In figure \ref{sle_segments} an SLE curve is shown for 
$\kappa=8/3$. There are two copies of the same curve. 
In the copy on the left it is divided into segments which 
correspond to equal changes in the capacity, while  
in the copy on the right it is divided into segments of 
equal fractal variation. 
When capacity is used to determine the segments, the segments appear
to have varying lengths, while the segments determined using the 
fractal variation appear to have equal length. 

In the next section we will consider this fractal variation in several
lattice models, and present Monte Carlo simulations that show
the limit in \reff{fvar} converges to a non-random constant.
In section 3 we use Monte Carlo simulations
to compare random variables defined using the natural parametrization 
in the lattice models with the corresponding random variables for 
SLE curves parametrized by their fractal variation. Section 4 gives 
some details about our simulations of SLE and the various lattice models.

\section{Fractal variation of discrete models}
\label{discrete}

We begin with our main conjecture:

\begin{conjecture}
Let $\gamma(t)$ be a random curve in the scaling limit of a critical 
lattice model with parametrization given by the length of the lattice
curve suitably scaled (the natural parametrization). 
Then the fractal variation of $\gamma[0,t]$ exists and 
is proportional to $t$. 
(The constant of proportionality will depend on the lattice.) 
\end{conjecture}

The fractal variation of $\gamma[0,t]$ is a priori 
a random variable. Part of the conjecture is that this random variable
is not random. Note that for models which are defined in the half plane,
the scaling limit is expected to be invariant under dilations. 
This implies that for models in the half plane, 
if the fractal variation is nonrandom, 
then it must be proportional to the natural parameterization.
In this section we will support the conjecture that the fractal variation 
is not random by numerically computing the fractal variation for the 
LERW, the SAW, Ising interfaces and percolation interfaces at a fixed 
value of their natural parameterization.
The simulations of the next section will test the conjecture that 
this fractal variation is proportional to the natural parameterization.
More details on the definition of these models
and our simulations may be found in section \ref{details}.

The fractal variation is the limit as $\Delta t \ra 0$ 
of $fvar(\gamma[0,t],\Delta t)$. The quantity $fvar(\gamma[0,t],\Delta t)$
is a random variable. For a particular lattice model we simulate this
random variable for several values of $\Delta t$, and for each 
value plot the cumulative distribution function, 
i.e., $P(fvar(\gamma[0,t]) \le x)$
as a function of $x$. If the fractal variation is indeed constant, then
these cumulative distribution functions should converge to a function 
that is $0$ left of the value $fvar(\gamma[0,t])$ and $1$ to the right. 

\begin{figure}[tbh]
\includegraphics{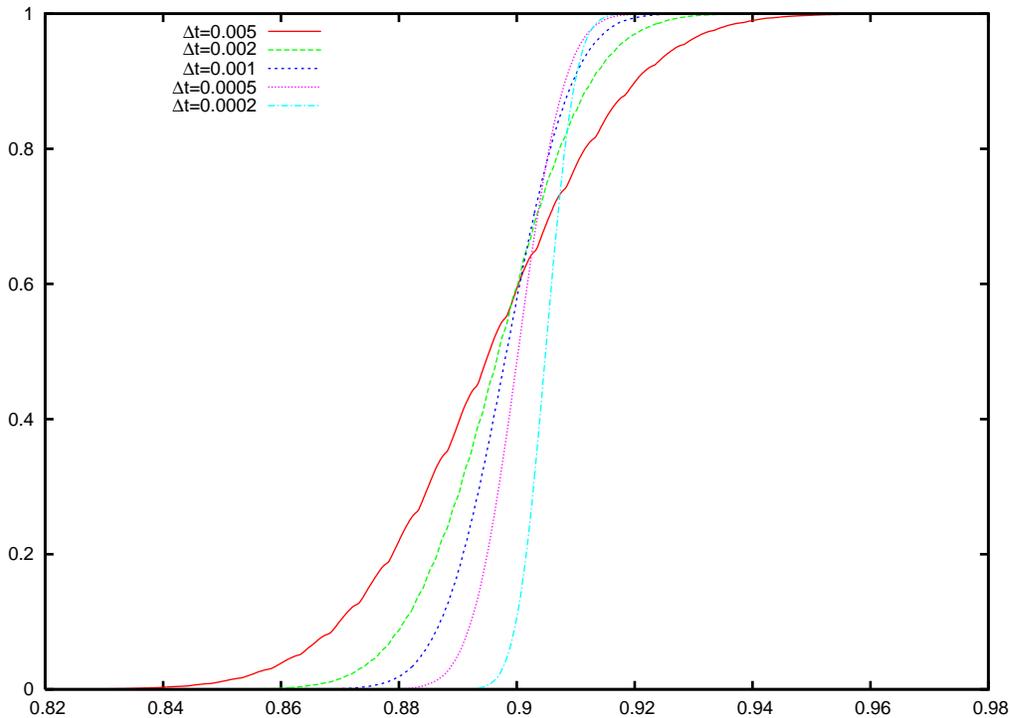}
\caption{\leftskip=25 pt \rightskip= 25 pt 
The cumulative distribution of the random variable 
$fvar(\gamma[0,t],\Delta t)$ for the LERW for a fixed $t$ and 
several values of $\Delta t$. (Color online.) 
} 
\label{lerw_pvar}
\end{figure}

Figure \ref{lerw_pvar} shows these cumulative distribution functions 
for the LERW
for several values of $\Delta t$. Note that the range of the horizontal 
axis shown is rather narrow and does not include zero. 
We should caution that we cannot take $\Delta t$ too small. It must be 
large enough that the number
of steps in the lattice walk corresponding to a single $\Delta t$ 
is large. If $\Delta t$ is extremely small, there
will be only a single step in the lattice walk corresponding to 
each $\Delta t$ and and  $fvar(\gamma[0,t],\Delta t)$  will 
just equal $n \Delta t$ where $n$ is the number of steps in 
the lattice walk. Thus as $\Delta t \ra 0$, we will first see 
$fvar(\gamma[0,t],\Delta t)$ converging to a constant, but then there will 
be a crossover where it begins to converge to a different constant. 
The beginning of this crossover is seen in figure \ref{lerw_pvar} 
where the curve with the smallest value of $\Delta t$, $2 \times 10^{-4}$, 
is shifted to the right of the other curves.

\begin{figure}[tbh]
\includegraphics{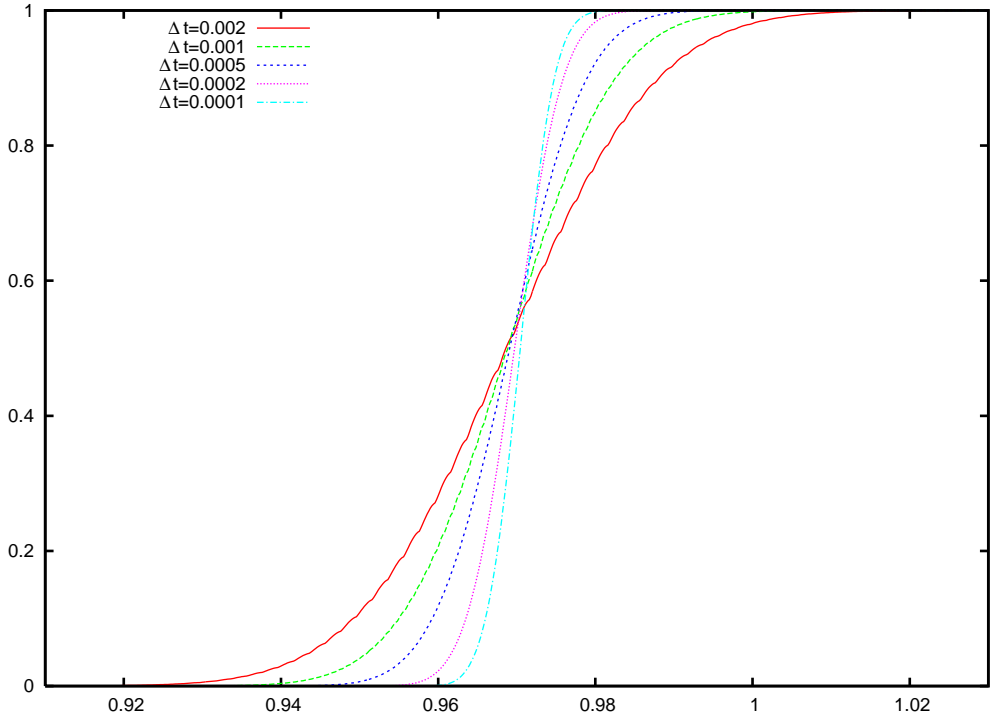}
\caption{\leftskip=25 pt \rightskip= 25 pt 
The cumulative distribution of the random variable 
$fvar(\gamma[0,t],\Delta t)$ for the SAW for a fixed $t$ and 
several values of $\Delta t$. (Color online.) 
} 
\label{saw_pvar}
\end{figure}

For the SAW the cumulative distributions of $fvar(\gamma[0,t],\Delta t)$ 
are shown in figure \ref{saw_pvar}. 
For the SAW it is possible to simulate
walks with a very large number of steps (one million). This is the reason
one does not see any shift of these curves even for the smallest $\Delta t$
of $10^{-4}$. 


\begin{figure}[tbh]
\includegraphics{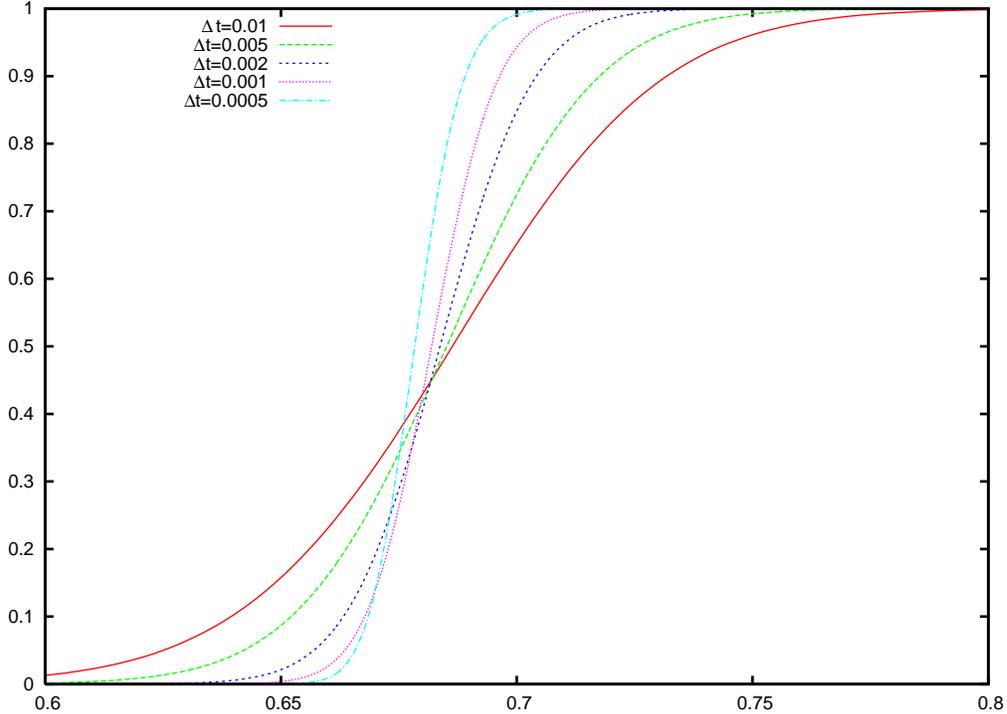}
\caption{\leftskip=25 pt \rightskip= 25 pt 
The cumulative distribution of the random variable 
$fvar(\gamma[0,t],\Delta t)$ for the Ising model for a fixed $t$ and 
several values of $\Delta t$. (Color online.) 
} 
\label{ising_pvar}
\end{figure}

For interfaces in the critical Ising model, figure \ref{ising_pvar}
shows the cumulative distributions of $fvar(\gamma[0,t],\Delta t)$. 
The smallest value
of $\Delta t$ shown is $5 \times 10^{-4}$, larger than the smallest value 
plotted for the LERW or SAW. One can already see the curve shifting to the 
left for this relatively large value of $\Delta t$. The Ising model is 
the most difficult to simulate since it is a truly two dimensional model.
So the interfaces we can generate are not nearly as long as for the 
other lattice models. 


\begin{figure}[tbh]
\includegraphics{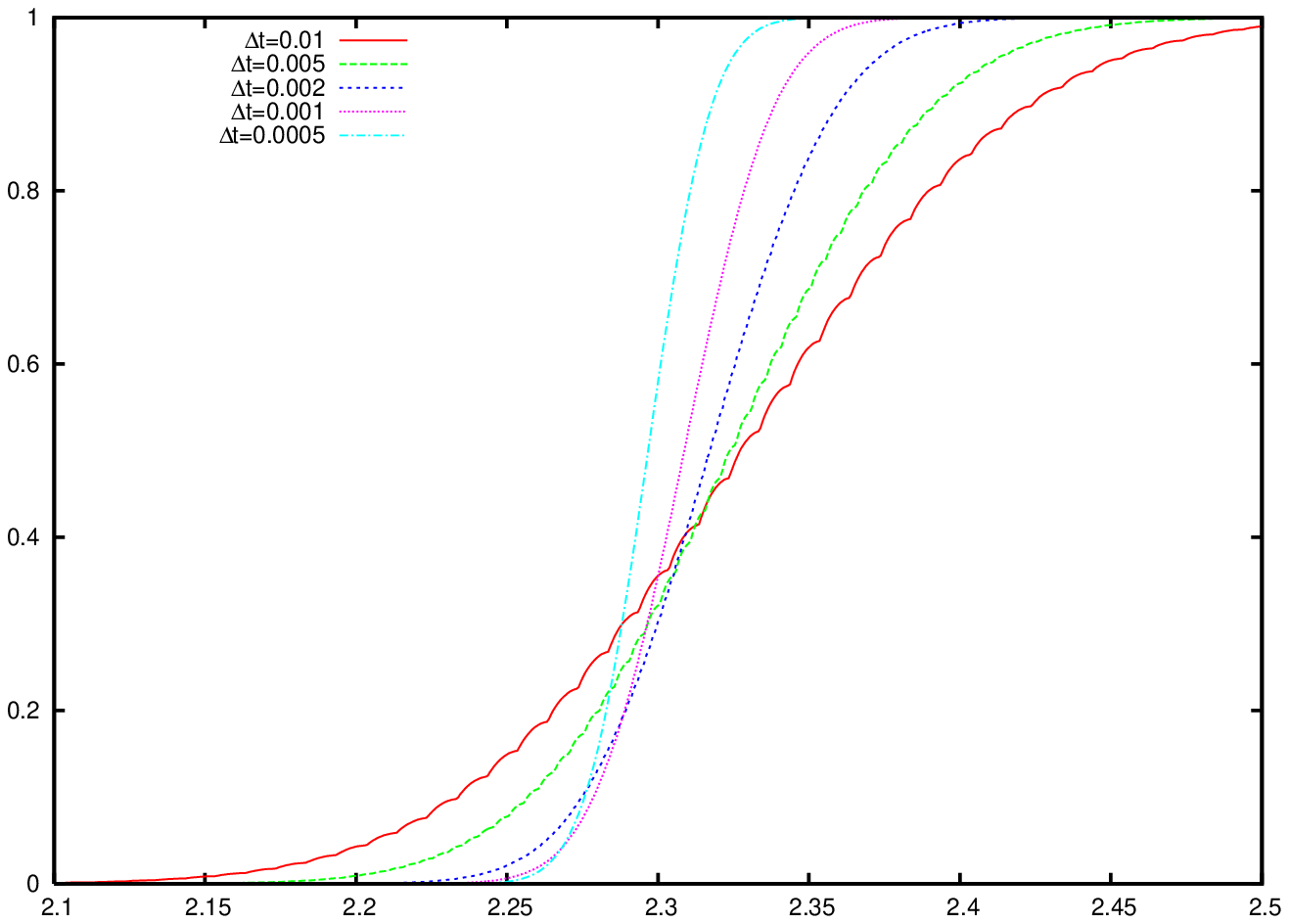}
\caption{\leftskip=25 pt \rightskip= 25 pt 
The cumulative distribution of the random variable 
$fvar(\gamma[0,t],\Delta t)$ for percolation for a fixed $t$ and 
several values of $\Delta t$. (Color online.) 
} 
\label{perc_pvar}
\end{figure}

The percolation interface is the easiest of the four lattice models to 
simulate, and we can generate curves with four million steps. However, 
this model has the largest Hausdorff dimension of the lattice models, and 
this makes the finite lattice effects occur at relatively large values
of $\Delta t$. 
Figure \ref{perc_pvar} shows the cumulative distributions of 
$fvar(\gamma[0,t],\Delta t)$.
For $\Delta t=5 \times 10^{-4}$, the smallest value shown, one can see the 
curve shifting to the left. In fact, this shift is just barely visible 
for $\Delta t=10^{-3}$.


\begin{figure}[tbh]
\includegraphics{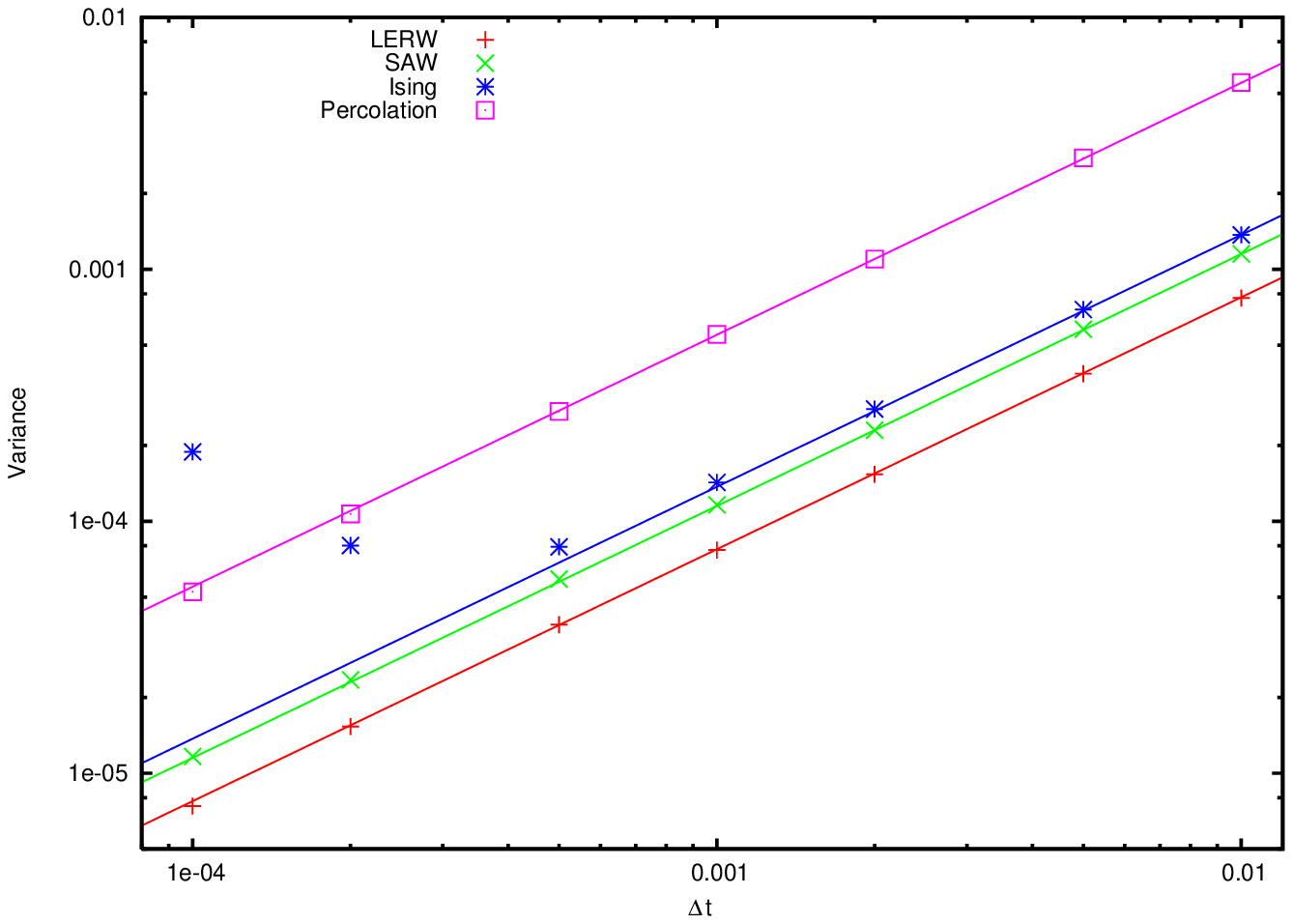}
\caption{\leftskip=25 pt \rightskip= 25 pt 
For each of the four lattice models we plot the variance of the 
random variable $fvar(\gamma[0,t],\Delta t)$ as a function of $\Delta t$. 
(Color online.) 
}
\label{pvar_variance}
\end{figure}

Figures \ref{lerw_pvar},\ref{saw_pvar},\ref{ising_pvar} and \ref{perc_pvar}
indicate that the random variables 
$fvar(\gamma[0,t],\Delta t)$ are converging to a constant as $\Delta t$
goes to zero. To study this quantitatively, we plot the variance of 
these random variables as a function of $\Delta t$ in figure 
\ref{pvar_variance}. We plot these variances for all four lattice models. 
In all four cases we also draw a line with slope 1 which attempts to fit the 
data for the values of $\Delta t$ before the crossover behavior sets in. 
This linear fit is extremely good until $\Delta t$ becomes small enough that 
the finite lattice effects are significant. 
Note that more values of $\Delta t$ are 
shown in figure \ref{pvar_variance} than were shown in the plots of the 
cumulative distributions of $fvar(\gamma[0,t],\Delta t)$.
The deviation of the data points in 
figure \ref{pvar_variance} from the corresponding straight line
begins roughly at the value of $\Delta t$ where one can see the 
plots of the cumulative distribution of 
$fvar(\gamma[0,t],\Delta t)$ beginning to shift in figures 
\ref{lerw_pvar},\ref{saw_pvar},\ref{ising_pvar} and \ref{perc_pvar}.

Since the data in figure \ref{pvar_variance} is well fit by a line 
with slope $1$, this indicates the variance of $fvar(\gamma[0,t],\Delta t)$ 
goes to zero as $\Delta t$. Note that the number of terms in the sum
defining $fvar(\gamma[0,t],\Delta t)$ is of the order of $1/\Delta t$. 
For the average of $N$ i.i.d. random variables, the variance also 
goes to zero as $1/N$. This suggests that the convergence of 
$fvar(\gamma[0,t],\Delta t)$ to a constant is some form of a law of 
large numbers. The random variables being summed are not independent, 
but one can hope that their correlations decay in some suitable way. 

\section{Comparison of SLE and the lattice models}
\label{sle_discrete}

The theorems and conjectures that state that the scaling limit of some
discrete model is SLE$_\kappa$ are usually statements that if we use the 
parametrization by capacity in both models, then the parametrized 
random curves in the two models have the same distribution. 
It is natural to expect that if we use the fractal variation to 
parametrize the curves in the two models, then these parametrized curves
should have the same distribution. Assuming that the fractal variation 
of the scaling limit of the discrete model is just a constant times the 
natural parametrization, if we use the natural parametrization 
for the discrete model (times a suitable constant) 
and the fractal parametrization for the SLE, then the parametrized curves 
should have the same distribution. We will test this conjecture
by considering random variables which depend on the parametrization of 
the curve. 

The most obvious random variable to study would be to consider the 
point on the random curve at a fixed value of the parametrization 
(the natural parametrization for the discrete model, 
the fractal variation parametrization for SLE). 
We conjecture these two parametrizations are equivalent, but there is 
a non-trivial constant of proportionality between them that we would
need to estimate. Instead, we will consider random variables that avoid
the need to compute this proportionality constant. 

In the following we use a superscript $'$ to indicate quantities
defined in lattice models, while quantities without such a superscript 
will indicate quantities in SLE. 
For example, $\gamma$ will denote an SLE curve while $\gamma^\prime$ denotes
a random curve from the scaling limit of a discrete model. 
$T$ will denote a value of the parameter in SLE where the parametrization
is defined using the fractal variation along the SLE curve, 
and $T^\prime$ will denote a value of the natural parametrization in 
a discrete model.

We consider three of the discrete models in the upper half plane:
the loop-erased random walk, the self-avoiding walk and the 
percolation exploration process. For these models we study the 
following random variables. Fix $R>0$. First we consider the SLE curve.
Look at the portion of the curve from its start at the origin until 
it first hits the semicircle of radius $R$. 
Let $T$ be the fractal variation of this part of the curve. 
(Of course, $T$ is random.) Let $\gamma(t)$ denote the curve parametrized
by the fractal variation. So $\gamma(T)$ is the first point on the  curve
with $|\gamma(T)|=R$. Consider the point $\gamma(T/2)$. It can be 
thought of as the point that is halfway along the curve from $0$ to 
$\gamma(T)$. 

For the discrete models, let $T^\prime$ be the time 
when it first hits the semicircle. (This is the number of steps 
suitably scaled.) Let $\gamma^\prime$ denote the curve
for the scaling limit of the discrete model and consider the 
point $\gamma^\prime(T^\prime/2)$. Then we expect that $\gamma(T/2)$ 
and $\gamma^\prime(T^\prime/2)$ have the same distribution.
Let $(X,Y)$ and $(X^\prime,Y^\prime)$ be these two random points. 
We test the conjecture by comparing the distributions of $X$ and $X^\prime$
and comparing the distributions of $Y$ and $Y^\prime$.
We also write the points in polar coordinates, $X+iY=R\exp(i\Theta)$, and
compare the distributions of $R$ with $R^\prime$ and of $\Theta$ with 
$\Theta^\prime$. 

\begin{figure}[tbh]
\includegraphics{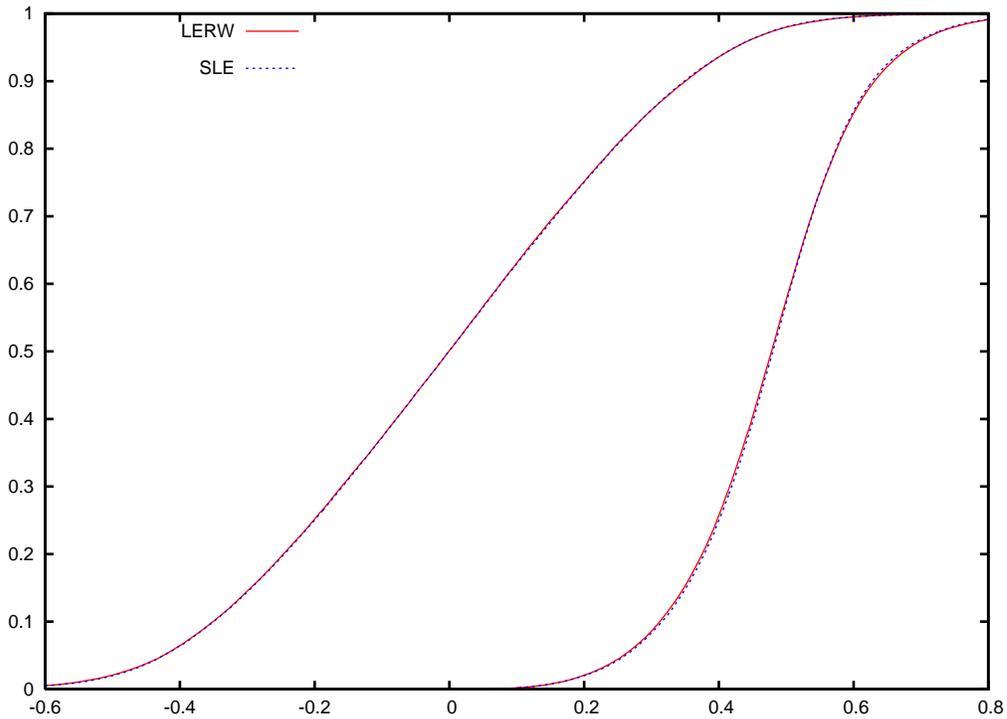}
\caption{\leftskip=25 pt \rightskip= 25 pt 
Cumulative distributions of $X$ coordinate (left two curves) and 
$Y$ coordinate (right two curves) of midpoint for LERW and SLE$_2$.
(Color online.) 
}
\label{lerw_sle_mid}
\end{figure}

\begin{figure}[tbh]
\includegraphics{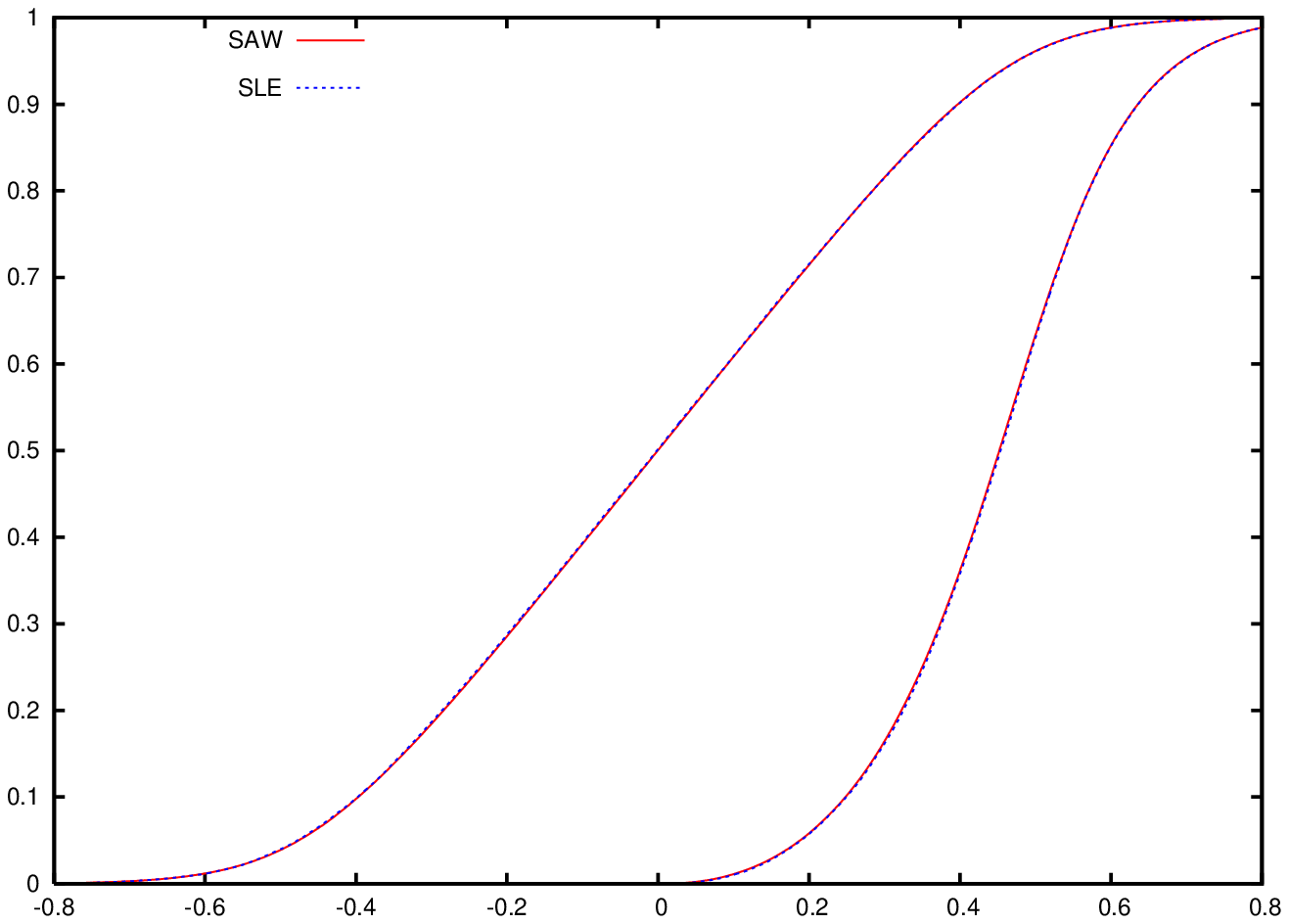}
\caption{\leftskip=25 pt \rightskip= 25 pt 
Cumulative distributions of $X$ coordinate (left two curves) and 
$Y$ coordinate (right two curves) of midpoint for SAW and SLE$_{8/3}$.
(Color online.) 
}
\label{saw_sle_mid}
\end{figure}

\begin{figure}[tbh]
\includegraphics{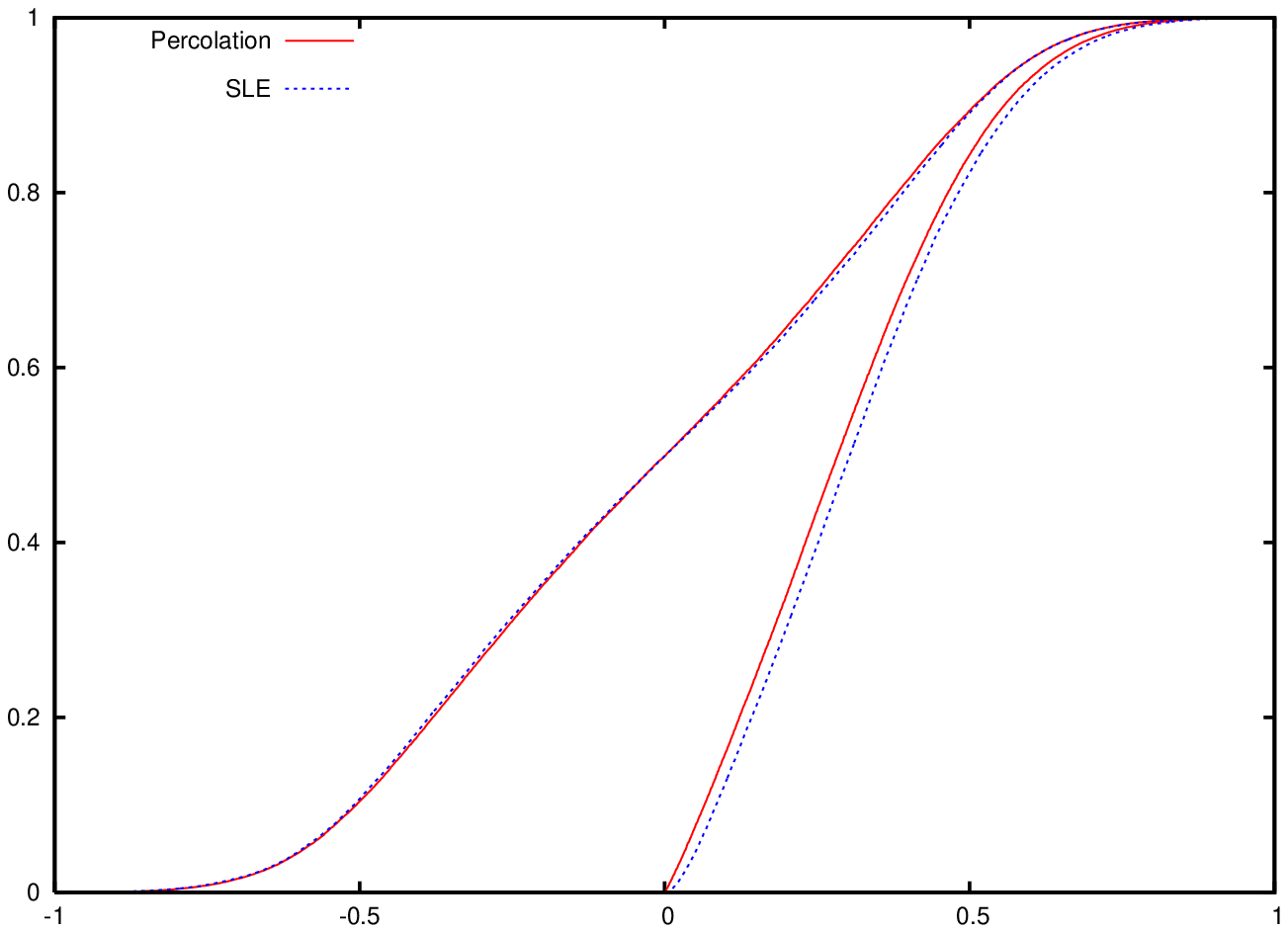}
\caption{\leftskip=25 pt \rightskip= 25 pt 
Cumulative distributions of $X$ coordinate (left two curves) and 
$Y$ coordinate (right two curves) of midpoint for percolation and SLE$_6$.
(Color online.) 
}
\label{perc_sle_mid}
\end{figure}

The cumulative distributions of $X,X^\prime,Y$ and $Y^\prime$ for the 
LERW, SAW and percolation are shown in figures \ref{lerw_sle_mid}, 
\ref{saw_sle_mid}, and \ref{perc_sle_mid}.
Each figure contains four curves. For the LERW and SAW, the 
cumulative distributions
of $X$ and $X^\prime$ agree so well that the two curves are virtually
indistinguishable in the figure. Likewise, the $Y$ and $Y^\prime$ curves
are virtually indistinguishable. So in figures  \ref{lerw_sle_mid} and
\ref{saw_sle_mid} it appears there are only two curves. For percolation
the difference between the cumulative distributions 
for percolation and SLE$_6$ 
are more noticeable, especially for $Y$ and $Y^\prime$. 
We do not fully understand this discrepancy, but believe it is caused 
by the difficulty of simulating SLE when $\kappa>4$. 

\begin{figure}[tbh]
\includegraphics{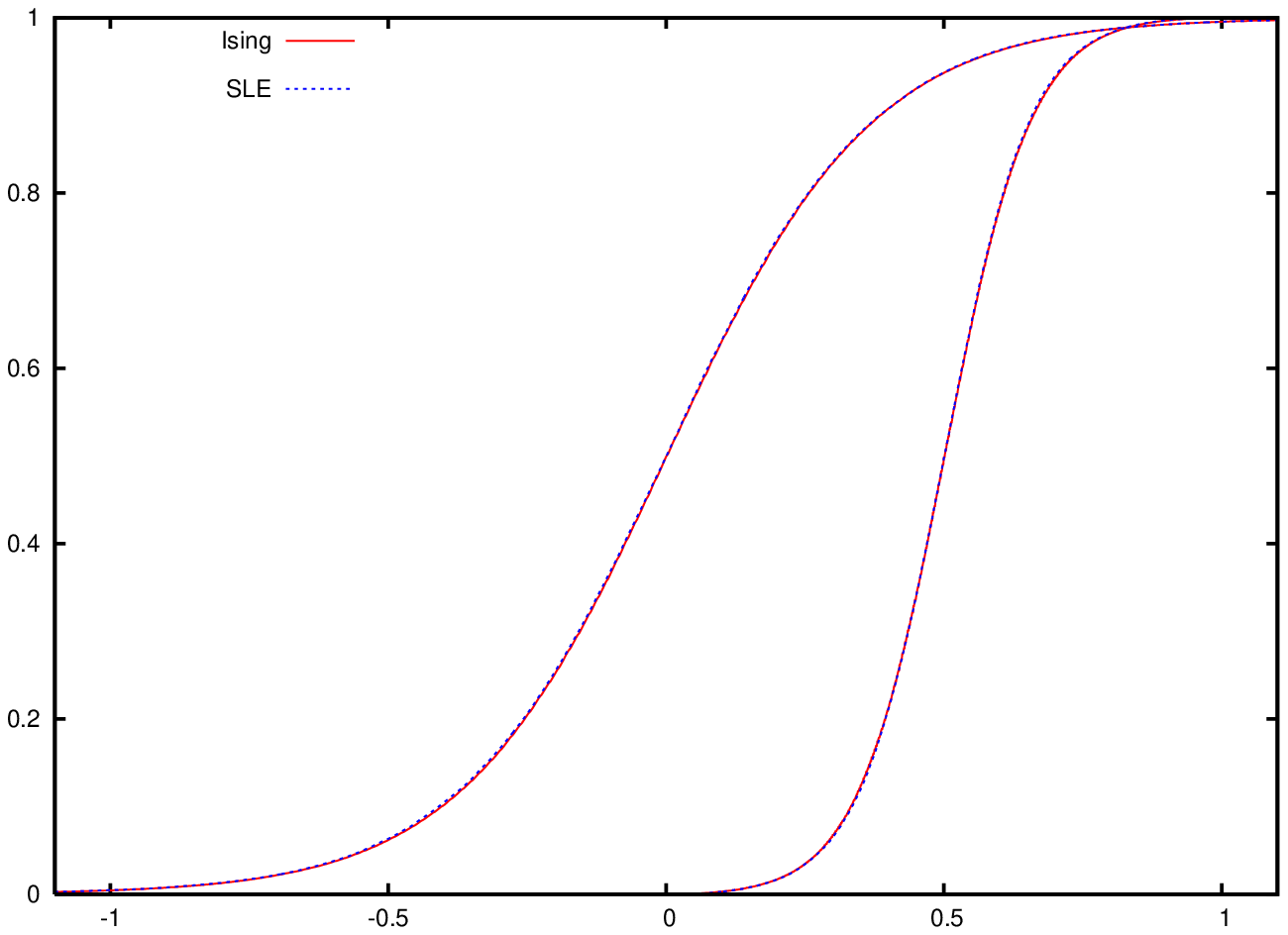}
\caption{\leftskip=25 pt \rightskip= 25 pt 
Cumulative distributions of $X$ coordinate (left two curves) and 
$Y$ coordinate (right two curves) of midpoint for Ising and SLE$_3$.
(Color online.) 
}
\label{ising_sle_mid}
\end{figure}


Table 1 shows the maximum difference between the 
cumulative distributions
for the discrete model and SLE. With two exceptions these differences 
are on the order of one percent or less. 
The two exceptions are the $Y,Y^\prime$ 
and $R$,$R^\prime$ comparisons for percolation. 
We do not show any plots of the cumulative distributions 
of $R,R^\prime, \Theta$
or $\Theta^\prime$. 

For the Ising model we consider a slightly different random variable 
since we cannot simulate the Ising model in a half plane. 
For the Ising model itself we take a rectangle of width $W$ and height $L$
with corners at $-W/2$, $W/2$, $-W/2+iL$ and $W/2+iL$ where 
$W$ is several times as large as $L$. 
On the top and bottom sides the boundary conditions are $-1$ 
for boundary sites with negative $x$ coordinate and  $+1$ for boundary 
sites with non-negative $x$ coordinate. We impose antiperiodic boundary
conditions between the vertical sides. These boundary conditions
force an interface from $0$ to $iL$ and approximate an infinite strip
\cite{bbh}. 
We let $T^\prime$ be the natural length of the Ising interface 
and consider $\gamma^\prime(T^\prime/2)$. 
We compare this Ising model with chordal SLE$_3$ 
in an infinite strip starting at $0$ and ending at $iL$, 
where $L$ is the width of the strip. 
We let $T$ be the fractal variation of the entire SLE curve
and consider the point $\gamma(T/2)$. 
Then as for the other models
we compare the cumulative distributions of the $x$ and $y$ coordinates of 
the random points $\gamma(T/2)$ and $\gamma^\prime(T^\prime/2)$.
Figure \ref{ising_sle_mid} compares the cumulative distributions 
of $X,X^\prime,Y$ and $Y^\prime$. 
Table 1 shows the maximum difference between their cumulative distributions,
and those of the polar coordinates. 

\begin{table}
\begin{center}
\begin{tabular}{|l|c|c|c|c|}
\hline
& X & Y & R & $\theta$ \\
\hline
LERW         &  0.003077  &    0.009211   &   0.007674   &  0.002425 \\ 
SAW          &  0.002306  &    0.005238   &   0.007831   &  0.002808 \\ 
Percolation  &  0.009293  &    0.040818   &   0.042395   &  0.013513 \\ 
Ising        &  0.003686  &    0.005446   &   0.003371   &  0.003141 \\ 
\hline
\end{tabular}
\caption{\protect \leftskip=25 pt \rightskip= 25 pt 
For a particular model and random variable, 
the table gives the maximum difference between the cumulative distribution
function for the random variable in the lattice model and in SLE with 
the corresponding value of $\kappa$. }
\end{center}
\end{table}

\section{Simulation details}
\label{details}

In this section we provide some details about the simulations of both 
the lattice models and SLE. 

To simulate SLE over the time interval $[0,t]$ we let 
$0=t_0 < t_1 < t_2 < \cdots < t_n=t$ be a partition of the time
interval. The times $t_k$  play a special role, but the random curves 
are still defined for all time. 
We replace the Brownian motion in the driving function by a stochastic
process that equals the Brownian motion at the 
times $t_k$, and is defined in between these times so that the Loewner 
equation may be solved explicitly. There are several ways to do this. 
We interpolate the driving function in between the time $t_{i-1}$ and 
$t_i$ by a square root function. 
If $g_t$ is the conformal map that takes the 
half plane minus the SLE hull back to the half plane, then  
this corresponds to approximating $g_t$ by a sequence of conformal maps,
each of which maps the half plane minus a linear slit back to the 
half plane. 

The simplest choice of the $t_i$ is to take a uniform partition of $[0,t]$.
This does not work well since it produces points along the SLE curve
that are far from being uniformly spaced. Instead we use an ``adaptive''
choice of the times \cite{rohde}.
Let $\Delta x>0$ and let $\gamma$ denote the SLE curve.
If $|\gamma(t_i)-\gamma(t_{i-1})| > \Delta x$ we divide the 
time interval into two equal halves. We repeat this process until 
$|\gamma(t_i)-\gamma(t_{i-1})| \le \Delta x$ for all the time intervals. 
Note that when we divide a time interval, we must use a Brownian bridge to 
choose the value of the random driving function at the time at the midpoint. 

Computing points along the SLE curve requires evaluating the composition
of approximately $N$ conformal maps. If we compute $N$ points along the 
curve, the total time needed for the computation will be $O(N^2)$. 
One can speed up the computation of the composition by approximating
the conformal maps by Laurent series. This leads to an algorithm 
that takes a time approximately $O(N^{1.4})$. This faster algorithm 
is explained in \cite{tk_sle}.

Our simulations of SLE for $\kappa=2,8/3$ and $6$ are done in the half 
plane.  We simulate each SLE curve until it hits the semicircle of radius 1. 
For $\kappa=2$, we set $\Delta x=0.002$ and generated $163,000$ samples.
For $\kappa=8/3$ we set $\Delta x=0.002$ and generated $132,000$ samples.
For $\kappa=6$ we set $\Delta x=0.005$ and generated $74,000$ samples.
The simulation for $\kappa=3$ is somewhat different since we need the SLE 
in a strip, not the half plane. 
We generate this SLE in the half plane, and apply 
the conformal map that takes the half plane to an infinite strip, 
sending the origin to the origin and $\infty$ to $i$. We compute the 
distance between consecutive points on the curve for comparison with 
$\Delta x$ after we have applied this conformal map.
We stop the simulation when the tip of the SLE is within a distance
$\Delta x$ of $i$. 
We set $\Delta x=0.005$ and generated $104,000$ samples. 

Lawler, Schramm and Werner proved that radial LERW converges to 
radial SLE$_2$. Radial refers to the fact that the random walk 
that one loop-erases 
begins at an interior point of the given domain and is conditioned to 
exit at a particular boundary point.
One can also consider the chordal LERW in which the random walk 
begins at a boundary point and is conditioned to remain in the domain 
and exit at a particular boundary point. Zhan showed this model 
converges to chordal SLE$_2$ \cite{zhan}. 

The LERW walk that we simulate is chordal LERW in the half plane from 
$0$ to $\infty$. This means that we take an ordinary random walk 
beginning at the origin and condition it to remain in the upper half 
plane. Then we erase the loops in chronological order. An ordinary 
random walk conditioned to remain in the upper half plane is known as 
the half plane excursion. It
is trivial to simulate since it is just given by a random walk
beginning at $0$ with transition probabilities that only depend on
the vertical component of the present location of the walk.
If the site has vertical component $k$, then the walk 
moves up with probability $(k+1)/4k$ and down with probability $(k-1)/4k$. 
The walk moves to the right or left with probability $1/4$. 
(See, for example, section 0.1 of \cite{lawler}).
The half-plane excursion is transient, i.e., each lattice site is visited 
by the excursion a finite number of times. This implies that the loop erasure
makes sense. (For a recurrent walk all parts of the walk would eventually
be part of a loop and so be erased.)
Note, however, that if we take an infinite half plane excursion and only 
consider the first $T$ steps and loop-erase this walk, the result will not 
agree with the loop-erasure of the full infinite excursion.
A site which is visited by the excursion before time $T$ may be erased 
by a loop formed after time $T$. 

In practice there is no way to know if a visit to a site will be erased 
by some future loop without simulating the entire excursion. 
So in the simulation we do the following. 
We generate a half-plane excursion, erasing the loops as 
they are formed. We stop when the resulting walk has $N$ steps. 
If $\alpha$ is small, then the distribution of our walk for the first
$\alpha N$ steps will be close to the true distribution of the 
first $\alpha N$ steps of the LERW.

For the simulations comparing the distribution of the ``midpoint'' 
with the corresponding point in SLE$_2$, we set $N=10^5$ and generated
118,000 samples.
We use a semicircle of radius $\rho N^{1/d_H}$ where 
$N=10^5$ is the number of steps. We set $\rho=0.1,0.2,0.4,0.8$.
If $\rho$ is too large, the finite length effects will begin to appear. 
We find no evidence of these effects for $\rho=0.1,0.2,0.4$, and good agreement
between these cases. The finite length effects are 
seen when $\rho=0.8$. The plots shown in figure \ref{lerw_sle_mid}
are for $\rho=0.4$. We rescale
all distances so that the radius of the semicircle is $1$. 
For the simulations of the fractal variation we do not see any finite
length effects and so we use the full $10^5$ steps in the walk in computing
the fractal variation. In these simulations we generated 132,000 samples. 

The SAW in the upper half plane is defined as follows. Let $N$ be a positive
integer. We consider all nearest neighbor walks with $N$ steps
in the upper half plane which begin at the origin and do not visit
any site more than once. Then we put the uniform probability measure on this
finite set of walks. We then let $N \rightarrow \infty$ to get a 
probability measure on infinite self-avoiding walks on the unit lattice
in the upper half plane. Finally, we take the lattice spacing to zero.  
Lawler, Schramm and Werner conjectured that this scaling limit is SLE$_{8/3}$ 
\cite{lsw_saw}. They proved the existence of the $N \rightarrow \infty$ 
limit, but the existence of the limit as the lattice spacing goes to zero
has not been established. Simulations of the SAW support their conjecture
\cite{tk_saw_sle_one,tk_saw_sle_two}.

The SAW in the half plane with a fixed number of steps may be simulated 
by the pivot algorithm, a Markov Chain Monte Carlo method. 
\cite{ms}. We use the fast implementation of this algorithm 
introduced in \cite{tk_pivot}. For the SAW there is an issue similar to
the LERW. The pivot algorithm produces the uniform distribution on 
the set of walks with $N$ steps. But this is not the distribution of 
the infinite SAW in the half plane restricted to walks of length $N$. 
As with the LERW, we address this problem by simulating walks with $N$ steps
but then working with random variables that typically only depend on a
relatively small initial part of the walk. 

For the simulations comparing the distribution of the ``midpoint'' 
with the corresponding point in SLE$_{8/3}$, we set $N=10^6$. 
We ran the pivot algorithms for $2 \times 10^9$ iterations of the 
Markov chain, and sampled the midpoint at each iteration. 
Of course, the resulting samples are far from independent. 
We use a semicircle of radius $\rho N^{1/d_H}$ with $\rho=0.2$. 
For the simulations for the fractal variation we again used walks 
with $N=10^6$ and ran the simulation for $2 \times 10^9$ iterations of the 
Markov chain, but we only computed the fractal variation every $10^3$ 
iterations since this computation is relatively time consuming. 
Again, we computed the fractal variation of the full $10^6$ step walk.
Earlier Monte Carlo studies of the fractal variation in the SAW and 
the comparison of the SAW and SLE$_{8/3}$ using the fractal variation
may be found in \cite{tk_saw_sle}

The Ising model we simulate is defined on a triangular lattice 
in a rectangle. Mixed boundary conditions are used on the 
top and bottom sides and antiperiodic boundary conditions
on the vertical sides as described in the previous section.
This forces an interface into the system that begins at the origin
and ends at $iL$. Smirnov has announced a proof that the scaling limit
of this interface is SLE$_3$ \cite{smirnov_icm,smirnov_ising}. 

We simulated the Ising model at its critical point with 
the Wolff algorithm \cite{wolff}.
The simulation was done in a rectangle of size $W=2000 \sqrt{3}$
by $L=1500$. We ran the Wolff algorithm for $1.4 \times 10^7$
iterations, sampling the random variable every 100 iterations.
The number of steps in the interface is random. We computed the 
fractal variation of the first $\rho L^{d_H}$ steps with 
$\rho=0.3$ and $0.6$. The data plotted is for $\rho=0.6$. 

The percolation model we study is site percolation on the triangular lattice
in the upper half plane, but we will describe it using the  
hexagonal lattice in the upper half plane. 
Each hexagon is colored white or black with probability $1/2$. 
The hexagons along the negative real axis are white and those along the 
positive real axis are black. This forces an interface which starts 
with the bond through the origin between the adjacent differently 
colored hexagons on the real axis. This interface is the unique
curve on the hexagon lattice which begins at this bond and 
has all white hexagons along one side of the interface and all black 
ones along the other side. Smirnov proved conformal invariance for 
this model, and so the interface is SLE$_6$ \cite{smirnov_perc}. 
See also Camia and Newman \cite{cn}. 

The percolation interface is the easiest of the lattice models to 
simulate. The key observation is that one should not generate 
the entire percolation configuration in a large rectangle 
but rather generate the configuration only as needed to determine the 
next step in the interface.
Interfaces with several million steps can be generated in a few seconds. 
Note that unlike the LERW or SAW there is no finite length effect.
If we generate interfaces with $n$ steps, they will have exactly the same 
distribution as the first $n$ steps of interfaces of length $N$ 
where $N>n$. 
The simulation was done with $N=4 \times 10^6$ 
and 98,000 samples were generated.
We used a semicircle of radius $\rho \times N^{1/d_H}$ with $\rho=0.2$.
For the simulations of the fractal variation, we also used
$N=4 \times 10^6$, but generated $206,000$ samples. 

\section{Conclusions}
\label{conclusions}

Our main conclusion is that for the various lattice models (LERW, SAW, 
Ising and percolation), the fractal variation exists and is proportional 
to the natural parametrization of the random curves defined using the 
length of the lattice curves. Our Monte Carlo
simulations support this conclusion in two ways. 
The first set of simulations show that the random variables 
$fvar(\gamma[0,t],\Delta t)$ converge to a constant as $\Delta t$
goes to zero. The second set of simulations support the conclusion in a 
more indirect way. If the conclusion is true for a lattice model and 
the scaling limit of that model is $SLE_\kappa$, then the fractal 
variation of the SLE curves should exist and provide a way to parametrize
the SLE curve that corresponds to the natural parametrization in the 
discrete model, up to a constant. Our second set of simulations tested this
by comparing the distribution of random variables that depend on the 
parametrization of the curve. 

\bigskip
\bigskip

\noindent {\bf Acknowledgments:}
Visits to the Kavli Institute for Theoretical Physics and the 
Banff International Research Station made possible many useful
interactions. The author thanks 
David Brydges, Greg Lawler, Don Marshall, Daniel Meyer, Yuval Peres, 
Stephen Rohde, Oded Schramm, Wendelin Werner and Peter Young for 
useful discussions. 
In particular, the definition of fractal variation we use grew out of 
these discussions.
This research was supported in part by the National Science Foundation 
under grants PHY99-07949 (KITP) and DMS-0501168 (TK).

\end{document}